\documentclass[12pt]{amsart}
\usepackage{amssymb}
\usepackage[mathscr]{eucal}
\usepackage{epsfig}
\usepackage{verbatim}
\usepackage{graphicx}
\usepackage{hyperref}

\newcommand{\C}{\ensuremath{\mathbb C}}
\newcommand{\R}{\ensuremath{\mathbb R}}
\newcommand{\PP}{\ensuremath{\mathbb P}}

\newcommand{\J}{\ensuremath{\mathcal J}}
\newcommand{\m}{\ensuremath{\mathfrak m}}

\newcommand{\Uaj}[1]{\ensuremath{U_{\alpha,#1}}}
\newcommand{\taj}[1]{t_{\alpha,#1}}
\newcommand{\uaj}[1]{u_{\alpha,#1}}
\newcommand{\vaj}[1]{v_{\alpha,#1}}

\DeclareMathOperator{\ord}{ord}

\DeclareMathOperator{\Bl}{Bl}

\DeclareMathOperator{\lct}{lct}

\theoremstyle{definition}
\newtheorem{defn}{Definition}[section]

\newtheorem{example}[defn]{Example}

\newtheorem{notation}[defn]{Notation}

\theoremstyle{plain}
\newtheorem{lem}[defn]{Lemma}
\newtheorem{induction}[defn]{Induction}
\newtheorem*{lem*}{Lemma}
\newtheorem{thm}[defn]{Theorem}

\newtheorem{prop}[defn]{Proposition}

\theoremstyle{remark}
\newtheorem{rem}[defn]{Remark}

\title{Multiplier ideals of general line arrangements in $\C^3$}
\author{Zachariah C.~Teitler}
\date{8/16/05}

\address{Department of Mathematics, SLU 10687, Hammond, LA 70402}
\email{zteitler@selu.edu}

\begin{document}

\bibliographystyle{plain}       

\begin{abstract}
We consider the multiplier ideals of the ideal of a reduced union
of lines through the origin in $\C^3$.
For general arrangements of lines, we calculate the multiplier ideals.
\end{abstract}

\maketitle

\section{Introduction}

Multiplier ideals have been used recently to answer several questions in algebraic
geometry (see, for example, \cite{siu}, \cite{el:theta}, or \cite{els:symbolic-powers}),
but they are hard to compute and relatively few examples are known
(see \cite{howald:monomial}, \cite{amanda},
\cite{blickle:mis-of-monomial-ideals-on-TVs},
\cite{mustata:hyperplane-arrangements},
and \cite{saito:mis-of-locally-conical-divs}).
In this paper we add to the list of known examples by computing multiplier ideals
in the case of a reduced union of lines through the origin in $\C^3$,
which we call an arrangement of lines.

Recall that for an ideal $I \subset \C[x_0,\dots,x_n]$,
regarded as an ideal on $\C^{n+1}$,
a \textbf{log resolution} of $I$ is a proper birational map
$f:X \to \C^{n+1}$, with $X$ smooth, such that the total transform
$I\cdot\mathscr{O}_X = \mathscr{O}_X(-F)$
is the ideal of a divisor $F$ with normal crossings support.
Then for $\lambda \geq 0$, $\lambda \in \R$, the $\lambda$th
\textbf{multiplier ideal}
$\J(I^\lambda)$ is given by
\[
  \J(I^\lambda)
    = f_{*} \mathscr{O}_X( K_{X/\C^{n+1}} - \lfloor \lambda F \rfloor) .
\]
More details on multiplier ideals may be found in \cite{pag2}.

Let $I$ be the ideal of an arrangement of lines
through the origin in $\C^3$.
In this case, the resolution process
begins naturally by blowing up the origin.
On the blowup space $X_0 = \Bl_0(\C^3)$, the pullback ideal
$I \cdot \mathscr{O}_{X_0}$ may have embedded components
supported in the exceptional divisor, making the resolution of $I$ non-trivial.
These embedded components correspond to certain subsets of $\PP^2$ containing
the set of points in $\PP^2$ which correspond to the lines in the arrangement.
We first consider this phenomenon.
For general arrangements, the embedded components are ``as nice as possible.''
There are three possibilities: there may be no embedded components,
there may be a smooth curve, or there may be a reduced set of points.
In each of these cases we describe the resolution and compute the multiplier ideal.

\section*{Acknowledgments}

I would like to thank my advisor, Rob Lazarsfeld,
for suggesting this question to me and for his guidance throughout
the preparation of my dissertation, of which this work is part.
I would also like to thank Tommaso de~Fernex
and Jessica Sidman for helpful comments and suggestions.

\section{The blowup of the vertex of an affine cone}

An arrangement of lines through the origin in $\C^3$ is
the affine cone over an arrangement of points in $\PP^2$,
and its log resolution begins by blowing up the origin.
As mentioned in the previous section, certain embedded components
may appear in the blowup space $\Bl_0(\C^3)$, depending on the geometry
of the arrangement of points.
In this section we will investigate this phenomenon in the slightly more
general setting of an arbitrary affine cone.

Let $Z \subset \PP^n$ be a nonempty closed subscheme
with saturated homogeneous ideal $I$.
Let $f:X = \Bl_0(\C^{n+1}) \to \C^{n+1}$ be the blowup of the origin,
with exceptional divisor $E$.
We wish to describe the pullback ideal $I \cdot \mathscr{O}_X$.
For this purpose we introduce some terminology.
We repeat some definitions and examples from
the companion paper to this one,~\cite{zct:smoothness-of-envelopes}.

\begin{defn}\label{defn:envelopes}
For $d\geq 0$, we define the \textbf{$d$th degree envelope},
or \textbf{$d$-envelope}, of $Z$ to be
the closed subscheme $Z_d \subset \PP^n$
defined by $I_d$, the degree $d$ piece of $I$.
The degree envelopes form a decreasing chain which begins with $\PP^n$
and stabilizes at $Z$.
If $Z_d \neq Z_{d-1}$, we say $d$ is a \textbf{geometric generating degree} of $I$.
\end{defn}

\begin{example}\begin{enumerate}
\item If $Z$ is a complete intersection of type $(d_1,\dots,d_r)$
with $d_1 < \dots < d_r$, then the geometric generating degrees of $I$
are exactly the $d_i$.
For each $i$, let $H_i$ be a hypersurface of degree $d_i$ such that
$Z = H_1 \cap \dots \cap H_r$.
Then $Z_{d_1} = H_1$, $Z_{d_2} = H_1 \cap H_2$, and so on.
\item Let $Z$ be five general reduced points in $\PP^2$.
Then $Z_2$ is the unique conic containing $Z$, and $Z_3 = Z$.
The geometric generating degrees are $2$ and $3$.
\item Let $Z$ be eight general reduced points in $\PP^2$.
Then there is a pencil of cubics passing through $Z$,
so $Z_3$ consists of the nine basepoints of this pencil.
That is, $Z_3$ is the union of $Z$ with an extra ninth point
(distinct from $Z$ because $Z$ is general).
The geometric generating degrees are $3$ and $4$.
\item Let $Z$ be four reduced points in $\PP^2$ with three collinear, but not all four.
Say the points $P_1$, $P_2$, and $P_3$ lie on the line $L$, and the point $P_4$
lies off of $L$.
Then $Z_2 = L \cup P_4$ and $Z_3 = Z$.
In this case a degree envelope has components of different dimensions.
The geometric generating degrees are $2$ and $3$.
\item Let $C$ be a smooth cubic and let $Z$ be eleven general reduced points on $C$.
Then $Z_3 = C$.
There is a unique point $P \in C$ such that $Z \cup P$ is the complete intersection
of $C$ with a quartic curve, and $Z_4 = Z \cup P$, twelve points ($P$ is distinct from
all the points of $Z$ by generality).
Finally, $Z_5 = Z$.
In this case, $I$ has three geometric generating degrees, $3$, $4$, and $5$.
\end{enumerate}
\end{example}

The following lemma will clarify the relationship between
the geometric generating degrees of $I$ and the usual
degrees of (algebraic) generators of $I$.

\begin{lem}\label{lem:ggds}
Let $I \subset S = \C[x_0,\dots,x_n]$ be a saturated homogeneous ideal.
Say $I = (H_1,\dots,H_s)$, with $H_i$ homogeneous, $\deg H_i = d_i$,
and $d_1 \leq \dots \leq d_s$.
Then:
\begin{enumerate}
\item every geometric generating degree of $I$ is one of the integers $d_i$.
\item $d_1$ is a geometric generating degree of $I$.\qed
\end{enumerate}
\end{lem}

\begin{rem}
We regard $\PP^n$ and $Z$ itself as trivial degree envelopes of $Z$
(We have $Z_d = \PP^n$ for $d < d_1$, in the notation of~\ref{lem:ggds},
and $Z_d = Z$ for $d \gg 0$.)
So $Z$ has no non-trivial degree envelopes if and only if $I$ has
only one geometric generating degree.
\end{rem}

As before let $I \subset S = \C[x_0,\dots,x_n]$ be  a saturated homogeneous ideal,
defining $Z \subset \PP^n$.
Consider the blowup of the origin, $f:X = \Bl_0(\C^n) \to \C^n$,
with exceptional divisor $E$.
We wish to relate the pullback ideal $I \cdot \mathscr{O}_X$
to the degree envelopes of $Z$.
We first introduce certain ideal sheaves on $X$.
For $d \geq 0$, let $I_d$ be the set of homogeneous forms of degree $d$ in $I$.
Note that $Z_d$ is cut out by $I_d$.
\begin{defn}\label{def:id}
For each $d\geq 0$, let $\mathscr{I}'_d$ be
the quotient ideal sheaf
\[
  \mathscr{I}'_d = (I_d \cdot \mathscr{O}_X :_{\mathscr{O}_X} \mathscr{O}_X(-d E)) .
\]
Also, let $\mathscr{I}'$ denote the ideal of the strict transform of
the affine cone over $Z$.
\end{defn}

Each $\mathscr{I}'_d$ is similar to, but not the same as,
the ideal of the strict transform of the
cone over $Z_d$. 
The $\mathscr{I}'_d$ satisfy the following properties.
\begin{lem}\label{lem:prop-of-id}
\begin{enumerate}
\item For each $d$, $I_d \cdot \mathscr{O}_X = \mathscr{O}_X(-d E) \cdot \mathscr{I}'_d$.
\item $\mathscr{I}'_d \subset \mathscr{I}'_{d+1}$.
\item $\mathscr{I}' = \sum_{d \geq 0} \mathscr{I}'_d = \mathscr{I}'_D$
for $D \gg 0$.
\item If $I_d S$ and $I_{d+1}S$ have the same saturation in $S$
(equivalently $Z_d = Z_{d+1} \subset \PP^n$),
then $\mathscr{I}'_d = \mathscr{I}'_{d+1}$.
\item If $I$ is generated by forms of degree $\leq D$,
then $\mathscr{I}' = \mathscr{I}'_D$.\qed
\end{enumerate}
\end{lem}

We can now give an expression for $I \cdot \mathscr{O}_X$.
At this point the expression includes redundant terms, but
we will see how to simplify it.

\begin{prop}\label{prop:blowup-of-vertex}
Let $I \subset S = \C[x_0,\dots,x_n]$ be a saturated homogeneous ideal,
let $f:X = \Bl_0(\C^{n+1}) \to \C^{n+1}$ with exceptional divisor $E$,
and let the ideal sheaves $\mathscr{I}'_d$ and $\mathscr{I}'$
be as in Definition~\ref{def:id}.
Let $I$ be generated (as an ideal)
by forms of degrees $d_1 < \dots < d_s$.
For $1 \leq j \leq s-1$, let
\[ \mathscr{P}_j = \mathscr{I}'_{d_j} + \mathscr{O}_X((d_1 - d_{j+1}) E) .\]
Then
\begin{equation}\label{eqn:primary-decomp-of-pullback}
  I \cdot \mathscr{O}_X
    = \mathscr{O}_X(-d_1 E) \cdot
        ( \mathscr{P}_1 \cap \dots \cap \mathscr{P}_{s-1} \cap \mathscr{I}' ) .
\end{equation}
\end{prop}

The following lemma will be useful.
\begin{lem}\label{lem:intersecting-ideals}
Let $R$ be a ring and $x$ a variable over $R$.
Let $0 \leq a_1 < a_2 < \dots < a_p$, and for each $1 \leq k \leq p$,
let $J_k$ be an ideal in $R$.
We have the following equality of ideals in $R[x]$:
\begin{multline*}
  x^{a_1} J_1^e + x^{a_2} J_2^e + \dots + x^{a_p} J_p^e \\
  = x^{a_1} \Big( (J_1^e + (x^{a_2-a_1})) \cap (J_1^e + J_2^e + (x^{a_3-a_1})) \cap \dots \\
    \quad \cap (J_1^e + \dots + J_{p-1}^e + (x^{a_p - a_1}))
    \cap (J_1^e + \dots + J_p^e) \Big) ,
\end{multline*}
where $J_k^e$ denotes the extended ideal $J_k R[x]$.\qed
\end{lem}

\begin{proof}[Proof of Proposition~\ref{prop:blowup-of-vertex}]
One can check $I \cdot \mathscr{O}_X \subset \mathscr{O}_X(-d_1 E) \cdot \mathscr{P}_j$
for each $j$ and $I \cdot \mathscr{O}_X \subset \mathscr{O}_X(-d_1 E) \cdot \mathscr{I}'$.
This gives the global inclusion $\subset$ in \eqref{eqn:primary-decomp-of-pullback}.

To show equality, it is enough to show that the global inclusion $\subset$
is an equality on the standard coordinate charts covering $X = \Bl_0(\C^{n+1})$.
Let $U_0 \subset X$ be the coordinate patch with coordinates $z_0,\dots,z_n$
such that $f:U_0 \to \C^{n+1}$ is given by
\[ (z_0,\dots,z_n) \mapsto (z_0,z_0z_1,\dots,z_0z_n) = (x_0,x_1,\dots,x_n) . \]
This corresponds to the map of rings
\begin{gather*}
  S = \C[x_0,\dots,x_n] \to S' := \C[z_0,\dots,z_n] , \\
  x_0 \mapsto z_0, \quad x_i \mapsto z_0 z_i, 1 \leq i \leq n .
\end{gather*}
We will show \eqref{eqn:primary-decomp-of-pullback} holds on $U_0$.
The other standard coordinate charts are similar.

We have $I = I_{d_1}S + \dots + I_{d_s}S$, so
\[ I \cdot \mathscr{O}_X(U_0) = I_{d_1}S' + \dots + I_{d_s}S' \]
is the ideal extension of $I$ in $S'$.
Let
\[
  I'_{d_i} = \{ \, H' = H(1,z_1,\dots,z_n) \mid H\in I_{d_i} \, \} \subset S' .
\]
Then $I_{d_i} S' = (z_0^{d_i}) I'_{d_i} S'$, so
\[
  I \cdot \mathscr{O}_X(U_0) = z_0^{d_1} I'_{d_1}S' + \dots + z_0^{d_s} I'_{d_s}S' .
\]
Similarly,
\[ \mathscr{P}_i(U_0) = I'_{d_i}S' + (z_0^{d_{i+1}-d_1}) , \]
and $\mathscr{I}'(U_0) = I'_{d_1}S' + \dots + I'_{d_s}S'$.

We apply Lemma~\ref{lem:intersecting-ideals},
with
$R = \C[z_1,\dots,z_n]$ and $x = z_0$, so $R[x] = S'$.
For $1 \leq i \leq p$, let $J_i$ be the ideal generated by $I'_{d_i}$ in $R$.
This completes the proof.
\end{proof}

Now we will apply Lemma~\ref{lem:ggds} to simplify our expression
for $I \cdot \mathscr{O}_X$.
It turns out we only have to include a term for each
\textit{geometric} generating degree of $I$.

\begin{prop}\label{prop:blowup-of-vertex-slimmer}
Let $I$ be a saturated homogeneous ideal in $S = \C[x_0,\dots,x_n]$,
let $X = \Bl_0(\C^{n+1})$ with exceptional divisor $E$,
and let the ideal sheaves $\mathscr{I}'_d$ and $\mathscr{I}'$
be as in Definition~\ref{def:id}.
Let $I$ have geometric generating degrees
$b_1 < \dots < b_q$, as in Definition~\ref{defn:envelopes}.
For $1 \leq j \leq q-1$, let
\[ \mathscr{Q}_j = \mathscr{I}'_{b_j} + \mathscr{O}_X((b_1 - b_{j+1}) E) .\]
Then
\begin{equation}\label{eqn:primary-decomp-of-pullback-slimmer}
  I \cdot \mathscr{O}_X
    = \mathscr{O}_X(-b_1 E) \cdot
        ( \mathscr{Q}_1 \cap \dots \cap \mathscr{Q}_{q-1} \cap \mathscr{I}' ) .
\end{equation}
\end{prop}

\begin{proof}
Let $d_1 < \dots < d_s$ be the degrees of all the generators
in some set of generators of $I$.
By Lemma~\ref{lem:ggds}, $b_1 = d_1$, and each $b_i$ is one
of the integers $d_j$.
For $1 \leq j \leq s-1$, let $\mathscr{P}_j$ be as in
Proposition~\ref{prop:blowup-of-vertex}.

If $d_j$ is not a geometric generating degree,
then $Z_{d_j} = Z_{d_{j-1}}$,
so $\mathscr{I}'_{d_j} = \mathscr{I}'_{d_{j-1}}$, by Lemma \ref{lem:prop-of-id}.
And since $d_j < d_{j+1}$, we have
\begin{equation}\label{eqn:redundant-pjs}
  \mathscr{P}_j =
  \mathscr{I}'_{d_j} + \mathscr{O}_X((d_1-d_{j+1})E)
  \subset
  \mathscr{I}'_{d_{j-1}} + \mathscr{O}_X((d_1-d_j)E) 
  = \mathscr{P}_{j-1} .
\end{equation}
Say $b_i = d_j < d_{j+1} < \dots < d_k < d_{k+1} = b_{i+1}$.
Then
\[ \mathscr{I}'_{b_i} = \mathscr{I}'_{d_j} = \dots = \mathscr{I}'_{d_k} . \]
Also, $\mathscr{O}_X((d_1-d_{k+1})E) = \mathscr{O}_X((b_1-b_{i+1})E)$.
This shows $\mathscr{Q}_i = \mathscr{P}_k$.
Chaining \eqref{eqn:redundant-pjs} together gives
\[
  \mathscr{Q}_i = \mathscr{P}_k \subset \mathscr{P}_{k-1}
    \subset \dots \subset \mathscr{P}_j .
\]
This shows that in \eqref{eqn:primary-decomp-of-pullback}
we can eliminate the redundant terms $\mathscr{P}_j,\dots,\mathscr{P}_{k-1}$,
leaving only $\mathscr{P}_k = \mathscr{Q}_i$.
\end{proof}

This description \eqref{eqn:primary-decomp-of-pullback-slimmer} of the pullback ideal
is not exactly a primary decomposition.
First of all, we have separated the divisorial and non-divisorial parts.
Second, the ideals $\mathscr{Q}_i$ are not necessarily primary.
In particular, $\mathscr{I}'$ corresponds to the strict transform of $Z \subset \C^n$,
so if $Z$ is reducible, then $\mathscr{I}'$ is not primary.
Similarly, if a degree envelope of $Z$ is reducible, the corresponding
ideal $\mathscr{Q}_i$ is not primary.
Finally, it may happen that some of these ideals are redundant.

It is still a useful description, and we may regard it as a first step
towards a primary decomposition of $I \cdot \mathscr{O}_X$.
We regard the $\mathscr{Q}_i$ as ``potential primary components''
of $I \cdot \mathscr{O}_X$.
They are supported along copies in $E \cong \PP^n$
of the degree envelopes of $Z$.
In this way, the degree envelopes of $Z$ can arise
as (potential) embedded components of the subscheme
defined by $I \cdot \mathscr{O}_X$.
These embedded components make the resolution of $I$ non-trivial.

\section{Degree envelopes of general sets of points in $\PP^2$}

The rest of this paper is devoted to finding a log resolution
and the multiplier ideals of an arrangement of lines in $\C^3$.
If the corresponding arrangement of points in $\PP^2$ has nontrivial
degree envelopes, then embedded components as in the previous section
will make the resolution non-trivial.
Therefore we only consider arrangements for which the degree envelopes
are ``as nice as possible.''
Let $Z$ be an arrangement of points in $\PP^2$ with ideal $I$.
We consider three cases:
\begin{itemize}
\item there is only one geometric generating degree of $I$,
\item there are two geometric generating degrees $d<e$,
and the intermediate degree envelope $Z_d$ is a smooth curve,
\item or there are two geometric generating degrees $d<e$,
and the intermediate degree envelope $Z_d$ is a set of reduced points.
\end{itemize}
In each case we will give the log resolution and multiplier ideals.
But before we do this, one might ask, which arrangements are covered by these
three cases?

One can show that for any $n\geq 1$, a general arrangement of $n$ lines in $\C^3$
falls into one of these three cases.
In the companion paper~\cite{zct:smoothness-of-envelopes}, we prove the following.
\begin{thm}\label{thm:ggd-of-general}
Let $n>1$.
Let $Z$ be a set of $n$ general points in $\PP^2$.
Let $d$ and $r$ be specified by $\binom{d+1}{2} \leq n = \binom{d+2}{2} - r$ with $r>0$,
so that $d$ is the lowest degree of a curve passing through $Z$,
and $r$ is the number of independent curves of degree $d$ passing through $Z$.
Let $I$ be the ideal of $Z$.
\begin{enumerate}
\item If $r=1$, the geometric generating degrees of $I$ are $\{d,d+1\}$
and the $d$-envelope $Z_d$ is a smooth curve of degree $d$.
\item If $r=2$ and $d>2$, the geometric generating degrees of $I$ are $\{d,d+1\}$,
and $Z_d$ is a set of $d^2$ distinct,
reduced points in $\PP^2$, a complete intersection of type $(d,d)$,
containing $Z$ together with $d^2-n = \binom{d-1}{2}$ extra points.
\item If $r=2$ and $d=2$ (so $n=4$),
then $2$ is the only
geometric generating degree of $I$.
\item If $r\geq 3$, then $d$ is the only
geometric generating degree of $I$.
\end{enumerate}
\end{thm}

In addition, these three cases include many ``special'' arrangements,
such as the transversal complete intersection of smooth curves of different degrees.

\section{Resolution and multiplier ideals of a line arrangement with one geometric generating degree}

We consider an arrangement $Z$ of lines through the origin of $\C^3$,
the corresponding set $\tilde{Z}$ of points in $\PP^2$,
and the saturated homogeneous ideal $I$.
Suppose $I$ has a single geometric generating degree.
This means the only degree envelopes of $\tilde{Z}$ are
just $\PP^2$ and $\tilde{Z}$ itself.
The resolution and multiplier ideals are essentially trivial in this case.

\begin{notation}
Let $\m$ denote the maximal ideal of the origin in $\C^3$.
\end{notation}

\begin{thm}\label{thm:no-emb-cpts}
We let $Z$, $\tilde{Z}$, and $I$ be as above.
If $I$ has only a single geometric generating degree $d$, then
the multiplier ideals of $I$ are as follows.
For $0 \leq \lambda < 2$,
$\J(I^\lambda) = \m^{\lfloor \lambda d \rfloor - 2}$.
For $2 \leq \lambda < 3$,
$\J(I^\lambda) = \m^{\lfloor \lambda d \rfloor - 2} \cap I$.
For $\lambda \geq 3$, we have the recursive description
\[ \J(I^\lambda) = I \J(I^{\lambda-1}) , \]
given by the Skoda theorem~\cite[Theorem~9.6.21]{pag2}.

In particular, the log canonical threshold of $I$ is
$\lct(I) = \min(3/d,2)$.
\end{thm}

\begin{proof}
Let $f_0:X_0 \to \C^3$ be the blowup of the origin, with exceptional
divisor $E_0$.
By hypothesis and Proposition~\ref{prop:blowup-of-vertex-slimmer},
$I \cdot \mathscr{O}_{X_0} = \mathscr{O}_{X_0}(-d E_0) \cdot \mathscr{I}'$,
where $\mathscr{I}'$ is the ideal sheaf of the strict transform $Z'$
of $Z$ in $X_0$.
Then $Z'$ is the disjoint union of strict transforms of the lines in $Z$.
For each line $\ell$ in $Z$, denote the corresponding component of $Z'$ by $\ell'$.
Each component of $Z'$ meets $E_0$ transversely.

\begin{rem}
In particular, $f_0:X_0 \to \C^3$ is a
\textbf{strong factorizing desingularization} of $Z$,
in the sense of Bravo and Villamayor~\cite{bravo-villamayor:strong}.
\end{rem}

Let $b:X \to X_0$ be the blowup of $Z'$.
Let the strict transform of $E_0$ in $X$ be denoted again by $E_0$.
Let the exceptional divisors of $b$ be denoted as follows:
for each line $\ell$ in $Z$,
let the irreducible exceptional divisor of $b$ over $\ell$
be denoted $E'_{\ell}$, and let $E' = \sum E'_{\ell}$.
Let $f = f_0 b : X \to \C^3$.
We see that $I \cdot \mathscr{O}_X = \mathscr{O}_{X}(-d E_0 - E')$
and $f$ is a log resolution of $I$.

The relative canonical divisor
$K_{X/\C^3} = \mathscr{O}_X( 2 E_0 + E' )$.
Therefore the $\lambda$th multiplier ideal of $I$ is
\[
\begin{split}
  \J(I^\lambda) &= f_*(\mathscr{O}_X( (2-\lfloor \lambda d\rfloor)E_0
                       + (1-\lfloor\lambda\rfloor) E' ) \\
                &= f_*(\mathscr{O}_X( (2-\lfloor \lambda d \rfloor)E_0)
                   \cap f_*(\mathscr{O}_X( (1-\lfloor\lambda\rfloor)E' ) \\
                &= m^{\lfloor \lambda d \rfloor - 2}
                   \cap I^{\langle \lfloor \lambda \rfloor - 1 \rangle} , \\
\end{split}
\]
where the second ideal is the $(\lfloor \lambda \rfloor -1 )$th symbolic power of $I$.
In particular, this is the unit ideal if $0 \leq \lambda < 2$
and $I$ if $2 \leq \lambda < 3$.
This proves the theorem.
\end{proof}

\section{Resolution and multiplier ideals of a line arrangement with a smooth curve envelope}

\begin{thm}\label{thm:emb-curve}
Let $\tilde{Z} \subset \PP^2$ be a finite set with saturated homogeneous ideal $I$.
Let $Z$ be the arrangement of lines corresponding to $\tilde{Z}$.
Suppose the geometric generating degrees of $I$ are $d,e$ with $d<e$
and the $d$-envelope $\tilde{Z}_d$ is a smooth curve.
Let $F_d$ be the homogeneous form of degree $d$ defining $\tilde{Z}_d$.
Then the multiplier ideals of $I$ are as follows.
For $0 \leq \lambda < 1$,
\[ \J(I^\lambda) = \m^{\lfloor \lambda d \rfloor - 2} . \]
For $1 \leq \lambda < 2$,
\[
  \J(I^\lambda)
    = \m^{\lfloor \lambda e \rfloor - (2 + e - d)}
    + \m^{\lfloor \lambda d \rfloor - (2 + d)} (F_d) .
\]
For $2 \leq \lambda < 3$,
\[
  \J(I^\lambda)
    = \big( \m^{\lfloor \lambda e \rfloor - (2+e-d)}
          + \m^{\lfloor \lambda e \rfloor - (2+2e-d)}(F_d)
          + \m^{\lfloor \lambda d \rfloor - (2+2d)}(F_d^2) \big)
      \cap I .
\]
For $\lambda \geq 3$, $\J(I^\lambda) = I \J(I^{\lambda-1})$,
by the Skoda theorem~\cite[Theorem~9.6.21]{pag2}.

In particular, the log canonical threshold of $I$ is
$\lct(I) = \min(3/d,(3+e-d)/e,2)$.
\end{thm}

\begin{proof}
Let $b_0:X_0 \to \C^3$ be the blowup of the origin, with exceptional
divisor $E_0$.
By the hypothesis and Proposition~\ref{prop:blowup-of-vertex-slimmer},
and using the notation $\mathscr{I}'$, $\mathscr{I}'_d$ as in that
proposition,
\[
  I \cdot \mathscr{O}_{X_0}
    = \mathscr{O}_{X_0}(-d E_0) \cdot
        \big( \mathscr{I}' \cap (\mathscr{I}'_d + \mathscr{O}_{X_0}((d-e)E_0)) \big).
\]

Let $C_d = \{F_d = 0\} \subset \C^3$.
Let $C_d^{(0)}$ be the strict transform of $C_d$ in $X_0$.
Then $C_d^{(0)}$ is smooth, and $C_d^{(0)}$ meets $E_0$
transversely along a smooth curve---a copy of $\tilde{Z}_d$ in $E_0 \cong \PP^2$.
The curve $E_0 \cap C_d^{(0)}$ is the support of an embedded component
of $I\cdot \mathscr{O}_{X_0}$.
We resolve this component first by blowing up repeatedly, before blowing up
the transforms of the lines in $Z$.

We construct a sequence of spaces obtained inductively by blowups.
We already have $b_0:X_0\to \C^3$. 
We may informally say $X_{-1} = \C^3$.
Suppose we have constructed $b_i:X_i \to X_{i-1}$, with exceptional divisor $E_i^{(i)}$,
and with $C_d^{(i)}$ the strict transform in $X_i$ of $C_d$.
Suppose also that $C_d^{(i)}$ is smooth and meets $E_i^{(i)}$ transversely
along a smooth curve.
Then let $b_{i+1}:X_{i+1} \to X_i$ be the blowup of this curve,
with exceptional divisor $E_{i+1}^{(i+1)}$.
Let $C_d^{(i+1)}$ be the strict transform in $X_{i+1}$ of $C_d$, or equivalently
of $C_d^{(i)}$.
Note that $C_d^{(i+1)}$ is smooth and meets $E_{i+1}^{(i+1)}$
transversely along a smooth curve.
We get a sequence
\[
  \cdots \to X_{i+1} \to X_i \to \cdots \to X_0 \to \C^3 .
\]
For $0 \leq j < i$ let $E_j^{(i)}$ be the strict transform
in $X_i$ of $E_j^{(j)}$.
Then $E_j^{(i)}$ meets only $E_{j \pm 1}^{(i)}$
and $C_d^{(i)}$ meets only $E_i^{(i)}$.

For each $i \geq 0$, the strict transform
$Z^{(i)}$ of $Z$ in $X_i$ consists of a union of pairwise disjoint
copies of lines, which each meet $E_i^{(i)}$
along the curve $E_i^{(i)} \cap C_d^{(i)}$,
and so they do not meet any $E_j^{(i)}$ with $j < i$.
Let $X = \Bl_{Z^{(e-d)}}(X_{e-d})$, with map $f:X \to \C^3$
given by the composition of all the blowdowns constructed previously.
For convenience, denote the strict transform of each $E_j^{(e-d)}$ in $X$ simply by $E_j$.
Let the irreducible exceptional divisor over each line $\ell$ of $Z$
be denoted $E'_{\ell}$, and let $\sum_{\ell} E'_{\ell} = E'$.

Then the $f$-exceptional locus $\sum_{j=0}^{e-d} E_j + E'$
is a simple normal crossings divisor.
We claim $f$ is a log resolution of the ideal $I$, with numerical data
\begin{equation}\label{eqn:numerical-data:emb-curve}
\begin{split}
  I \cdot \mathscr{O}_X = \mathscr{O}_X( -d E_0 -(d+1)E_1 - \dots - e E_{e-d} - E') , \\
  K_{X/\C^3} = 2 E_0 + 3 E_1 + \dots + (e-d+2) E_{e-d} + E' .
\end{split}
\end{equation}
To show this, we go by induction on $j$, with the following statement.
\begin{induction}\label{induction}
For $0 \leq j \leq e-d$, the space $X_j$ 
is covered by open sets $\{ \Uaj{j} \}_{\alpha}$ with local coordinates
$(\taj{j},\uaj{j},\vaj{j})$ satisfying the following conditions.
\begin{enumerate}
\item $E_j^{(j)} \cap \Uaj{j}$ is cut out by $\taj{j} = 0$ ;
\item $C_d^{(j)} \cap \Uaj{j}$ is cut out by $\uaj{j} = 0$ ;
\item The strict transform $Z^{(j)}$ either does not meet $\Uaj{j}$,
or there is a single line $\ell$ in $Z$ whose strict transform $\ell^{(j)}$
meets $\Uaj{j}$, and is cut out in that neighborhood by
$\uaj{j} = \vaj{j} = 0$ ;
\item
\[
  I \cdot \mathscr{O}_{X_j}|_{\Uaj{j}} =
  \begin{cases}
    (\taj{j}^{d+j})(\uaj{j},\taj{j}^{e-d-j}) ,      & \Uaj{j} \cap Z^{(j)} = \emptyset \\
    (\taj{j}^{d+j})(\uaj{j},\vaj{j} \taj{j}^{e-d-j}) , & \Uaj{j} \cap Z^{(j)} \neq \emptyset
  \end{cases}
\]
\item $F_d |_{\Uaj{j}} = \taj{j}^{d+j} \uaj{j}$ ;
\item For $H \in S$ homogeneous and $F_d\nmid H$, we have
$H |_{\Uaj{j}} = \taj{j}^{\deg H} h_j$, where
$h_j \notin (\taj{j},\uaj{j})$ .
\end{enumerate}
\end{induction}
\begin{rem}
For $j > e-d$, the statement still holds, replacing negative powers $\taj{j}^{e-d-j}$
with $1$, but we will not use this.
\end{rem}
\begin{proof}[Proof of induction]
For $j=0$, the existence of such coordinates follows
from the Implicit Function Theorem~\cite{gh},
since the hypersurfaces $E_i^{(i)}$ and $C_d^{(i)}$
meet transversally.
The expression for $I \cdot \mathscr{O}_{X_0}$ in local coordinates
follows from Proposition \ref{prop:blowup-of-vertex-slimmer}.
The expression for the pullback of $F_d$ is clear.
The expression for the pullback of $H$ can be proved as follows.
Since $H$ is homogeneous, $h|_{E_0^{(0)} \cap \Uaj{0}}$
is the same as $H$ on (an open subset of) $\PP^2$.
Then $H \notin (F_d)$ just means $h \notin ( \taj{0},\uaj{0} )$.
This proves the initial case $j=0$ of the induction.

The map $b_{j+1}:X_{j+1} \to X_j$ is, over $\Uaj{j}$,
the blowup of the line $\taj{j}=\uaj{j}=0$.
The preimage of $\Uaj{j}$ in $X_{j+1}$ is covered with two coordinate
charts in the usual way for blowups of linear subspaces.
One of these charts meets the exceptional divisors $E_j^{(j+1)}$
and $E_{j+1}^{(j+1)}$, but not $C_d^{(j+1)}$.
On this chart, the pullback of $I$
and the pullbacks of $F_d$ and $H$ are given by
the information in the induction statement.
We have to examine the other chart, which we denote $\Uaj{j+1}$,
with the same indexing set $\{\alpha\}$,
and show that these charts satisfy the induction hypothesis.

They clearly cover the curve $E_{j+1}^{(j+1)} \cap C_d^{(j+1)}$.
The chart $\Uaj{j+1}$ has coordinates $(\taj{j+1},\uaj{j+1},\vaj{j+1})$
such that the blowdown map $\Uaj{j+1} \to \Uaj{j}$ is given by
\[ (\taj{j+1},\uaj{j+1},\vaj{j+1}) \mapsto (\taj{j+1},\taj{j+1} \uaj{j+1},\vaj{j+1}) . \]
Then $E_{j+1}^{(j)}$, $C_d^{(j+1)}$, and $Z^{(j+1)}$ (if it meets $\Uaj{j+1}$)
are clearly cut out by the coordinates on $\Uaj{j+1}$ as claimed.
The pullback of $I$ from $\Uaj{j}$ to $\Uaj{j+1}$
gives the claimed generators.
Again, the pullback of $F_d$ clearly has the claimed description.
The pullback of $H$ also has the claimed description,
and since $h_j$ contains a term of the form $c \vaj{j}^a$ for some nonzero
constant $c$ and $a \geq 0$, we see that $h_{j+1}$ contains the term
$c \vaj{j+1}^a$, so $h_{j+1} \notin (\taj{j+1},\uaj{j+1})$, as claimed.

This completes the proof of the induction.
\end{proof}

Given this induction, we claim
\[
  I \cdot \mathscr{O}_{X_{e-d}}
    = \mathscr{O}_{X_{e-d}}(-d E_1 - \dots - e E_{e-d} ) \cdot \mathscr{I}' .
\]
(Here $\mathscr{I}'$ is the ideal sheaf of the proper transform $Z^{(e-d)}$.)
Away from $C_d^{(e-d)}$, we use the inductive description of
the pullback $I \cdot \mathscr{O}_{X_j}$
on the neighborhood $\Uaj{j}$ to see what $I \cdot \mathscr{O}_{X_j}$ is near $E_j$.
In a neighborhood of $E_{e-d}^{(e-d)} \cap C_d^{(e-d)}$, we get the claim
from the inductive description of $\mathscr{I}$ on $\Uaj{e-d}$.
And away from the exceptional locus, the claim is tautologically true.

\begin{rem}
This shows $f_{e-d}:X_{e-d} \to \C^3$ is a
\textbf{strong factorizing desingularization} of $Z$,
in the sense of Bravo and Villamayor~\cite{bravo-villamayor:strong}.
\end{rem}

Blowing up the lines $\ell' = \ell^{(e-d)}$ in $X_{e-d}$ then
gives a log resolution, as claimed.
We get exactly the numerical data for $I$ in
\eqref{eqn:numerical-data:emb-curve}.

We have shown that
\[
  \J(I^\lambda)
    = f_* \mathscr{O}_X \Big( \sum_{j=0}^{e-d} ((2+j)-\lfloor \lambda (d+j) \rfloor) E_j
            + ( 1 - \lfloor \lambda \rfloor ) E' \Big) .
\]
The ideal $I$ is homogeneous,
the map $f$ is $\C^{*}$-equivariant
(where $\C^{*}$ acts on $\C^3$ by scaling),
and each exceptional divisor $E_j$, $E'_{\ell}$ is stable
under the lifted $\C^{*}$ action on $X$
(see~\cite[Theorem 3.1]{equivariant-resolutions}).
Therefore the multiplier ideal $\J(I^\lambda)$ is also homogeneous.
So, to determine $\J(I^\lambda)$, it is sufficient
to characterize its homogeneous elements.
We do this by considering the valuations induced on $S=\C[x,y,z]$ by the 
irreducible $f$-exceptional divisors.

Naturally the valuation on $S$ induced by $E_0$ measures the
order of vanishing at the origin in $\C^3$,
and the valuation induced by each $E'_{\ell}$ is the order of vanishing
along the line $\ell$.
For the $E_j$ with $1 \leq j \leq e-d$,
the inductive statement~\ref{induction} shows that the valuation
determined by $E_j$ is as follows.
For homogeneous $H \in S$, $F_d \nmid H$, and $0 \leq j \leq e-d$,
we have $\ord_{E_j}(H) = \deg H$.
For $0 \leq j \leq e-d$, $\ord_{E_j}(F_d) = d + j$.

Every homogeneous element of $S$ has a unique factorization $H F_d^a$
with $F_d \nmid H$.
Hence $H F_d^a \in \J(I^\lambda)$ if and only if
$\deg H$, $a$ satisfy the inequalities
\begin{equation}\label{eqn:inequalities:emb-curve}
  \deg H + (d+j)a \geq \lfloor \lambda (d+j) \rfloor - (2+j) ,
  \quad 0 \leq j \leq e-d ,
\end{equation}
and
\[
  H F_d^a \in f_* \mathscr{O}_X((1-\lfloor \lambda \rfloor) E')
    = I^{\langle \lfloor \lambda \rfloor - 1 \rangle},
\]
the $(\lfloor \lambda \rfloor - 1)$th symbolic power of $I$.
In particular this is the unit ideal $(1)$ for $0 \leq \lambda < 2$,
and $I$ for $2 \leq \lambda < 3$.

Ignoring for the moment this second condition, we analyze the
inequalities \eqref{eqn:inequalities:emb-curve}
by fixing $d$, $a$, and $\lambda$:
then we get the single inequality
\[
\begin{split}
  \deg H
    & \geq \max_{0 \leq j \leq e-d} \lfloor \lambda (d+j) \rfloor - (2+ad+(a+1)j) \\
    & =
    \begin{cases}
      \lfloor \lambda d \rfloor - (2 + ad) , & \lambda \leq a+1 \\
      \lfloor \lambda e \rfloor - (2 + (a+1)e - d) , & \lambda \geq a+1
    \end{cases}
\end{split}
\]
Now determining which homogeneous forms $H F_d^a$ lie in the multiplier ideal
$\J(I^\lambda)$ is fairly mechanical.

Let $0 \leq \lambda < 1$.
Then $H F_d^a \in \J(I^\lambda)$ if and only if
\eqref{eqn:inequalities:emb-curve} holds for $\deg H$ and $a$.
Since $\lambda \leq a+1$ is automatic,
\eqref{eqn:inequalities:emb-curve} holds
if and only if $\deg H F_d^a \geq \lfloor \lambda d \rfloor - 2$.
Therefore $\J(I^\lambda) = \m^{\lfloor \lambda d \rfloor - 2}$.

Let $1 \leq \lambda < 2$.
Then $\deg H$, $a$ satisfy \eqref{eqn:inequalities:emb-curve}
if and only if $a=0$ and $\deg H  \geq \lfloor \lambda e \rfloor - (2+e-d)$
or $a \geq 1$ and $\deg H \geq \lfloor \lambda d \rfloor - (2+ad)$.
Therefore $H F_d^a \in \J(I^\lambda)$ if and only if
$H F_d^a = H \in \m^{\lfloor \lambda e \rfloor - (2+e-d)}$
or $a \geq 1$ and
$H F_d^a \in \m^{\lfloor \lambda d \rfloor - (2+ad)}(F_d^a) \subset \m^{\lfloor \lambda d \rfloor - (2+d)}(F_d)$.
Therefore
$\J(I^\lambda) = \m^{\lfloor \lambda e \rfloor - (2+e-d)} + \m^{\lfloor \lambda d \rfloor - (2+d)}(F_d)$.

Let $2 \leq \lambda < 3$.
In this case, the symbolic power $I^{\langle 1 \rangle} = I$.
Then $H F_d^a \in \J(I^\lambda)$ if and only if $\deg H$, $a$
satisfy \eqref{eqn:inequalities:emb-curve} and $H F_d^a \in I$.
Since $F_d \in I$, the condition $H F_d^a$ only matters if $a=0$.
We see that $\deg H$, $a$ satisfy \eqref{eqn:inequalities:emb-curve}
if and only if: $a=0$ and $\deg H  \geq \lfloor \lambda e \rfloor - (2+e-d)$,
$a = 1$ and $\deg H \geq \lfloor \lambda e \rfloor - (2+2e+d)$,
or $a \geq 2$ and $\deg H \geq \lfloor \lambda d\rfloor - (2+ad)$.
One checks that this corresponds to the ideal given in the statement
of the theorem.

This completes the proof of the theorem.
\end{proof}

\section{Resolution and multiplier ideals of a line arrangement with a set of points envelope}

\begin{thm}\label{thm:emb-pts}
Let $\tilde{Z} \subset \PP^2$ be a finite set with homogeneous ideal $I$.
Suppose the geometric generating degrees of $I$ are $d,e$ with $d<e$
and the $d$-envelope $\tilde{Z}_d$ is finite and reduced,
$\tilde{Z} \subset \tilde{Z}_d$.
Let $\widetilde{W} = \tilde{Z}_d \setminus \tilde{Z}$.
Let $I_W$ be the saturated ideal defining $\widetilde{W}$.
The multiplier ideals of $I$ are as follows.
For $0 \leq \lambda < 2$,
\[ \J(I^\lambda) = \m^{\lfloor \lambda d \rfloor - 2} . \]
For $2 \leq \lambda < 3$,
\[
  \J(I^\lambda)
    = ( \m^{\lfloor \lambda d \rfloor - 2} \cap I_W
        + \m^{\lfloor \lambda e \rfloor - 2(1+e-d)} ) \cap I .
\]
For $\lambda \geq 3$, $\J(I^\lambda) = I \J(I^{\lambda-1})$,
by the Skoda theorem~\cite[Theorem~9.6.21]{pag2}.

In particular, the log canonical threshold of $I$ is
$\lct(I) = \min(3/d,2)$.
\end{thm}

\begin{proof}
The proof is similar to the proof of Theorem~\ref{thm:emb-curve}.
On blowing up the origin in $\C^3$, we get embedded points.
These have to be resolved by blowing up $e-d$ times,
each time at points in the most recent exceptional divisor.
One checks similarly that this gives a log resolution of $I$.
The computation of pushforwards is also similar.
\end{proof}

\section{Examples}

\begin{example}
Let $Y \subset \PP^2$ be a set of $3$ non-collinear points.
Then by Theorem~\ref{thm:ggd-of-general}
the ideal of $Y$ has a single geometric generating degree, which is $2$.
By Theorem~\ref{thm:no-emb-cpts}, the line arrangement which is the affine
cone over $Y$ has log canonical threshold equal to $3/2$.

On the other hand, let $Z \subset \PP^2$ be a set of $3$ collinear points.
Then the ideal of $Z$ has the geometric generating degrees $1$ and $2$,
and $Z_1$ is the line on which $Z$ lies.
By Theorem~\ref{thm:emb-curve}, the line arrangement which is the affine
cone over $Z$ has log canonical threshold equal to $5/3$.

Therefore, even though $Z$ is more special than $Y$, the corresponding
line arrangement is quantitatively less singular.
This does not contradict the semicontinuity of log canonical thresholds,
since the degeneration is not flat.
On deforming three independent lines to
three coplanar lines, the flat limit consists of the three coplanar lines
together with an embedded point at the origin~\cite[pg.~72]{eisenbud-harris}.
\end{example}

\begin{example}
Similarly, let $Y \subset \PP^2$ be a set of $6$ points in general position.
By Theorem~\ref{thm:ggd-of-general} and Theorem~\ref{thm:no-emb-cpts},
the line arrangement which is the affine
cone over $Y$ has log canonical threshold equal to $1$.

Let $Z \subset \PP^2$ be a set of $6$ points on a smooth conic.
By Theorem~\ref{thm:emb-curve}, the line arrangement which is the affine
cone over $Z$ has log canonical threshold equal to $4/3$.

Again, even though $Z$ is more special than $Y$, the corresponding
line arrangement has a higher log canonical threshold.
As before, the deformation is not flat.
\end{example}

\bibliography{../biblio}   

\begin{thebibliography}{10}

\bibitem{blickle:mis-of-monomial-ideals-on-TVs}
Manuel Blickle.
\newblock Multiplier ideals and modules on toric varieties.
\newblock {\em Math. Z.}, 248(1):113--121, 2004.

\bibitem{bravo-villamayor:strong}
A.~Bravo and O.~Villamayor~U.
\newblock A strengthening of resolution of singularities in characteristic
  zero.
\newblock {\em Proc. London Math. Soc. (3)}, 86(2):327--357, 2003.

\bibitem{el:theta}
Lawrence Ein and Robert Lazarsfeld.
\newblock Singularities of theta divisors and the birational geometry of
  irregular varieties.
\newblock {\em J. Amer. Math. Soc.}, 10(1):243--258, 1997.

\bibitem{els:symbolic-powers}
Lawrence Ein, Robert Lazarsfeld, and Karen~E. Smith.
\newblock Uniform bounds and symbolic powers on smooth varieties.
\newblock {\em Invent. Math.}, 144(2):241--252, 2001.

\bibitem{eisenbud-harris}
David Eisenbud and Joe Harris.
\newblock {\em The geometry of schemes}, volume 197 of {\em Graduate Texts in
  Mathematics}.
\newblock Springer-Verlag, New York, 2000.

\bibitem{equivariant-resolutions}
Santiago Encinas and Orlando Villamayor.
\newblock A course on constructive desingularization and equivariance.
\newblock In {\em Resolution of singularities (Obergurgl, 1997)}, volume 181 of
  {\em Progr. Math.}, pages 147--227. Birkh\"auser, Basel, 2000.

\bibitem{gh}
Phillip Griffiths and Joseph Harris.
\newblock {\em Principles of algebraic geometry}.
\newblock Wiley Classics Library. John Wiley \& Sons Inc., New York, 1994.
\newblock Reprint of the 1978 original.

\bibitem{howald:monomial}
J.~A. Howald.
\newblock Multiplier ideals of monomial ideals.
\newblock {\em Trans. Amer. Math. Soc.}, 353(7):2665--2671 (electronic), 2001.

\bibitem{amanda}
Amanda~A. Johnson.
\newblock {\em Multiplier ideals of determinantal ideals}.
\newblock PhD thesis, U. Michigan, 2003.

\bibitem{pag2}
Robert Lazarsfeld.
\newblock {\em Positivity in algebraic geometry. {II}}, volume~49 of {\em
  Ergebnisse der Mathematik und ihrer Grenzgebiete. 3. Folge. A Series of
  Modern Surveys in Mathematics [Results in Mathematics and Related Areas. 3rd
  Series. A Series of Modern Surveys in Mathematics]}.
\newblock Springer-Verlag, Berlin, 2004.
\newblock Positivity for vector bundles, and multiplier ideals.

\bibitem{mustata:hyperplane-arrangements}
Mircea Musta{\c t}{\v a}.
\newblock Multiplier ideals of hyperplane arrangements.
\newblock To appear in TAMS., November 2004.

\bibitem{saito:mis-of-locally-conical-divs}
Morihiko Saito.
\newblock Multiplier ideals, b-function, and spectrum.
\newblock arXiv:math.AG/0402363, April 2004.

\bibitem{siu}
Yum-Tong Siu.
\newblock Invariance of plurigenera.
\newblock {\em Invent. Math.}, 134(3):661--673, 1998.

\bibitem{zct:smoothness-of-envelopes}
Zachariah~C. Teitler.
\newblock On the intersection of the curves through a set of points in
  $\mathbb{P}^2$.
\newblock arXiv:math.AG/0508307, 2005.

\end{thebibliography}

\end{document}